\documentclass{article}
\usepackage{amsmath}
\usepackage{graphicx}
\usepackage{amssymb}
\usepackage{amsthm}
\usepackage[english]{babel}
\usepackage{pgf}
\usepackage{tikz}
\usetikzlibrary{intersections}

\newcommand{\lrbrace}[1]{\left \lbrace #1 \right \rbrace }
\newcommand{\lrrb}[1]{\left ( #1 \right ) }
\newcommand{\skp}[1]{ \left \langle #1 \right \rangle }
\newcommand{\abs}[1]{\left| #1 \right|}

\newcommand{\te}[1]{\text{#1}}

\newcommand{\set}[1]{\mathbb{#1}}

\newtheorem{theorem}{Theorem}[section]
\newtheorem{definition}[theorem]{Definition}

\newtheorem{proposition}[theorem]{Proposition}

\begin{document}

\author{Sebastian Helmensdorfer}

\title{Space-time constructions for the mean curvature flow in a Ricci flow background}

\maketitle

\begin{abstract}

Given a solution of the (backwards) Ricci flow one can construct a so called canonical soliton metric on space-time (see \cite{estherpeter1, estherpeter2}). We observe that for a mean curvature flow within a (backwards) Ricci flow background, the space-time track of the mean curvature flow yields a canonical soliton of the coupled flow within the canonical Ricci soliton. We show that this provides a link between the differential Harnack estimate for the mean curvature flow and the differential Harnack estimate for the Ricci flow (see \cite{hamiltonrf, hamiltonmcf}). Moreover the second fundamental form of our canonical soliton matches the boundary term of the evolution of J. Lott's modified $\mathcal{F}$-functional for a Ricci flow with boundary (see \cite{lottmcfrfbckgr}). This functional also appears in quantum gravity (see \cite{yorkfunc, gibbonshawkingfunc}).

\end{abstract}

\section{Introduction}

Canonical solitons for the Ricci flow have been introduced by E. Cabezas-Rivas and P. Topping (see \cite{estherpeter1, estherpeter2}), building up on results by B. Chow, S. Chu and D. Knopf (see \cite{chowchu1, chowchu, chowknopf}). These are space-time constructions, which give geometric insight into optimal transportation in relation to Ricci flow (see \cite{lottoptimal, mccanntopping, toppingoptimal}) and into Harnack estimates for the Ricci flow (see \cite{hamiltonrf, brendleharnack, estherpeter2}). The idea can be adapted to curvature flows in Euclidean space, which provides geometric insight into the corresponding Harnack estimates (see \cite{kotschwarharnack}).

In this work we show that the construction of canonical solitons extends in a natural way to the case of mean curvature flow in a Ricci flow background. Given a solution of the (backwards) Ricci flow one can construct the canonical shrinking, steady (for the backwards Ricci flow) and expanding Ricci solitons (for the Ricci flow), which are metrics on space-time. We show that the space-time track of the associated mean curvature flow is itself a canonical soliton within one of the canonical Ricci flow solitons. 

This provides a relation between Hamilton's Harnack quantity for the mean curvature flow and Hamilton's Harnack quantity for the Ricci flow (see \cite{hamiltonrf, hamiltonmcf}). The first can be seen as the second fundamental form of our canonical soliton, whereas the latter corresponds to the Ricci curvature of the canonical expanding Ricci soliton (see \cite{estherpeter2, kotschwarharnack}).

Moreover J. Lott has introduced a version of G. Perelman's $\mathcal{F}$-functional (see \cite{perelman1}) for the Ricci flow with boundary (see \cite{lottmcfrfbckgr}), which is a weighted version of the Gibbons-Hawking-York functional from quantum gravity (see \cite{gibbonshawkingfunc, yorkfunc}). We show that the evolution of this $\mathcal{F}$-functional, which has also close relations to the Harnack quantities, corresponds to the second fundamental form of our canonical soliton.

\section{Preliminaries}

We want to study solutions of the mean curvature flow in a Ricci flow background. Consider families $F_t : M^{n} \rightarrow O^{n+1}$, $t \in I$ ($I$ a real interval), of smooth immersions in a Riemannian manifold $O^{n+1}$ equipped with a corresponding family of metrics $g^O (t)$ such that

\begin{eqnarray}
  \label{eq:meancurvatureflowinaricciflowbackground}
   \frac{\partial}{\partial t} F_t = - H \nu \\
   \nonumber 
   \frac{\partial}{\partial t} g^O = S
\end{eqnarray}

where 

$$S = - 2 Ric_{ g^O } \te{ (Ricci flow background) } $$
or 
$$S = 2 Ric_{ g^O } \te{ (backwards Ricci flow background) }.$$
 
We denote $M_t = F_t \lrrb{ M^n }$, $\nu$ denotes a choice of unit normal vector on $M_t$ and $H$ the mean curvature of $M_t$. The induced metric and the second fundamental form of $M_t$ are denoted by $g$ and $h$. We use coordinates $\lrrb{x^i}_{i = 1, \ldots, n}$ on $M^n$ and coordinates $\lrrb{y^\alpha}_{\alpha = 1, \ldots, n + 1}$ on $O^{n+1}$. For the curvature on $O^{n+1}$ we use the sign convention

\begin{equation}
  \nonumber
  \mathcal{ R }_{g^O} \lrrb{X,Y} = {}^{g^O} \nabla_Y {}^{g^O} \nabla_X - {}^{g^O} \nabla_X {}^{g^O} \nabla_Y + {}^{g^O} \nabla_{[X,Y]}.
\end{equation}

If $S = 2Ric$ we use $\tau$ instead of $t$ as the reverse time parameter, in line with standard notation.

\section{Canonical solitons}

To define the notion of a soliton in our setting, we follow \cite{lottmcfrfbckgr}. 

Let $\lrrb{O^{n+1}, g^O (t), f^O}$ for $t \in I$, $f^O : O^{n+1} \times I \rightarrow \set{ R }$ be a gradient Ricci flow soliton, i.e.

\begin{eqnarray}
	\label{eq:definitionmeancurvatureflowinaricciflowbackground}
	\frac{\partial}{\partial t} g^O = - 2Ric_{ g^O } \\
	\nonumber Ric_{ g^O } + Hess_{g^O} \lrrb{ f } + \frac{c}{2t} g^O = 0 \\
	\nonumber \frac{\partial f^O}{\partial t} = \abs{ {}^{g^O} \nabla f^O}^2
\end{eqnarray}

where $c = 0$ in the steady case ($I = \set{R}$), $c = 1$ in the expanding case ($I = (0,\infty)$) and $c = -1$ in the shrinking case ($I = (-\infty,0)$). 

\begin{definition}
 \label{def:mcfsolitonsinarfbackground}

 At a given time $t$ an immersed hypersurface $M_t$ in $O^{n+1}$ is a soliton of the mean curvature flow (within the Ricci soliton) if

 \begin{equation}
   \label{eq:mcfsolitoninarfbackground}
   H \pm \nu f^O = 0.
 \end{equation}

\end{definition}

Definition \ref{def:mcfsolitonsinarfbackground} can be understood as follows. Suppose $\lrrb{O^{n+1}, g^O (t), f^O}$ is a gradient steady Ricci soliton. Let $\lrrb{\psi_t}_{t \in I}$ be the one-parameter family of diffeomorphisms generated by the time-independent vector field ${}^{g^O} \nabla f^O$ and with $\psi_0 = Id$, such that $g^O (t) = \psi_t^* g^O (0)$. 

Then if a hypersurface $M_0$ satisfies (\ref{eq:mcfsolitoninarfbackground}) with a ``$+$'' sign, the hypersurfaces $M_t = \psi_t \lrrb{ M_0 }$ provide a solution of (\ref{eq:meancurvatureflowinaricciflowbackground}) up to tangential diffeomorphisms. Changing the sign in (\ref{eq:mcfsolitoninarfbackground}) means that the solution $M_t$ moves along the vector field $- {}^{g^O} \nabla f^O$.

For an expanding or shrinking Ricci soliton background we choose the one-parameter of diffeomorphisms generated by $\pm \frac{1}{\tau (t)} {}^{g^O} \nabla f^O$ for $\tau(t) = c t + 1 > 0$ and we get an analogous statement as above.

Examples for solutions of (\ref{eq:mcfsolitoninarfbackground}) and more generally examples for solutions of (\ref{eq:meancurvatureflowinaricciflowbackground}) with a Ricci soliton background can be found in \cite{tsatismcfinarfbck}.

\subsection{Canonical solitons within expanding Ricci flow solitons}

We consider solutions $\lrrb{M^n, M_t, O^{n+1}, g^O (t)}$ of (\ref{eq:meancurvatureflowinaricciflowbackground}) for $S = -2Ric_{g^O}$ and we choose $I = (0,T]$.

The construction of canonical expanding solitons can be described as follows (see \cite[theorem 2.1]{estherpeter2}). Suppose $\lrrb{O^{n+1}, g^O (t)}$, $t \in (0,T]$, has uniformly bounded curvature. For a large enough parameter $N > 0$ we define the metric $\check{g}^O$ on $O^{n+1} \times (0,T]$ by

\begin{equation}
 \label{eq:canonicalricciexpanders}
 \check{g}^O_{ij} = \begin{cases}
                     \frac{g^O_{ij}}{t}, \; i,j \geq 0 \\
		     0, \; i \geq 1, j = 0 \\
		     \frac{N}{2t^3} + \frac{R_{ g^O }}{t} + \frac{n + 1}{2 t^2}
                    \end{cases}
\end{equation}

where $0$ denotes the time direction on $O^{n+1} \times (0, T]$ and $R_{g^O}$ denotes the scalar curvature.

Then up to errors of order $\frac{1}{N}$ the metric $\check{g}$ is an approximate gradient expanding Ricci soliton:

\begin{equation}
 \nonumber
 E_N = Ric_{\check{g}} + Hess_{\check{g}} \lrrb{- \frac{N}{2t}} + \frac{1}{2} \check{g} \approx 0
\end{equation}

by which we mean that the quantity $N \abs{ E_N }_{\check{g}}$ is bounded uniformly locally on $O^{n+1} \times (0,T]$ (independently of $N$).

We denote the potential of the canonical expanding Ricci soliton by $\check{ f } \lrrb{y , t} = - \frac{N}{2t}$. We can show that there is a natural corresponding construction for the coupled flow.

\begin{theorem}
 \label{thm:canonicalmcfsolitonsinrfexpanders}

 Suppose $\lrrb{M^n, M_t, O^{n+1}, g^O (t)}$ is a solution of (\ref{eq:meancurvatureflowinaricciflowbackground}) such that $g^O(t)$ has uniformly bounded curvature and $M_t$ has uniformly bounded second fundamental form.

Define $\check{ \Sigma } = \lrbrace{ \lrrb{ x, t }: x \in M_t, t \in (0,T] }$ to be the space-time track of $M_t$ within the canonical Ricci expander $\lrrb{O^{n+1} \times (0, T], \check{g}}$. Then $\check{ \Sigma }$ is an approximate soliton, i.e.

\begin{equation}
 \nonumber
 \check{ E }_N = H^{ \check{ \Sigma } }  - \nu^{ \check{ \Sigma } }   \check{f} \approx 0
\end{equation}

where $H^{ \check{ \Sigma } }$ and $\nu^{ \check{ \Sigma } }$ denote the mean curvatue and the normal vector of $\check{ \Sigma }$. By this we mean that $N \abs{\check{ E }_N}_{\check{g}}$ is bounded uniformly locally on $\check{ \Sigma }$ (independently of $N$).

\end{theorem}
 
The proof of theorem \ref{thm:canonicalmcfsolitonsinrfexpanders} consists of computing the relevant quantities (see appendix \ref{appendix:canonicalmcfsolitonsexpandingrfbck}). We keep the terminology and call $\check{ \Sigma }$ a canonical soliton of the mean curvature flow (within a canonical expanding Ricci soliton).

\subsection{Canonical solitons within shrinking Ricci flow solitons}

We consider solutions $\lrrb{M^n, M_\tau, O^{n+1}, g^O (\tau)}$ of (\ref{eq:meancurvatureflowinaricciflowbackground}) for $S = 2Ric_{g^O}$ and we choose $I = (0,T]$.

Associated to the backwards Ricci flow $g^O (\tau)$ we have a canonical shrinking Ricci soliton which can be described as follows (see \cite[theorem 1.1]{estherpeter1}). Suppose $\lrrb{O^{n+1}, g^O (\tau)}$, $\tau \in (0,T]$, has uniformly bounded curvature. For a large enough parameter $N > 0$ we define the metric $\hat{g}^O$ on $O^{n+1} \times (0, T]$ by

\begin{equation}
 \label{eq:canonicalriccishrinkers}
 \hat{g}_{ij} = \begin{cases}
                     \frac{g^O_{ij}}{\tau}, \; i,j \geq 0 \\
		     0, \; i \geq 1, j = 0 \\
		     \frac{N}{2\tau^3} + \frac{R_{ g^O }}{\tau} - \frac{n + 1}{2 \tau^2}
                    \end{cases}
\end{equation}

where $0$ denotes the time direction on $O^{n+1} \times (0, T]$ and $R_{g^O}$ denotes the scalar curvature.

Then up to errors of order $\frac{1}{N}$ the metric $\hat{g}$ is an approximate gradient shrinking Ricci soliton:

\begin{equation}
 \nonumber
 E_N = Ric_{\hat{g}} + Hess_{\hat{g}} \lrrb{ \frac{N}{2\tau}} - \frac{1}{2} \check{g} \approx 0
\end{equation}

by which we mean that the quantity $N \abs{ E_N }_{\hat{g}}$ is bounded uniformly locally on $O^{n+1} \times (0, T]$ (independently of $N$). It has been shown in \cite{estherpeter1}, that the metric $\hat{g}$ relates to $\mathcal{L}$-optimal transportation (for details see \cite{estherpeter1, toppingoptimal}). 

We denote the potential of the canonical shrinking Ricci soliton by $\hat{ f } \lrrb{y , \tau} = \frac{N}{2\tau}$. As in the expanding case we can show that there is a natural corresponding construction for the coupled flow. The following result was first observed in a slightly different setting in \cite[section 6]{estherpeter1}. 

\begin{theorem}
 \label{thm:canonicalmcfsolitonsinrfshrinkers}

 Suppose $\lrrb{M^n, M_\tau, O^{n+1}, g^O (\tau)}$ is a solution of (\ref{eq:meancurvatureflowinaricciflowbackground}) such that $g^O(\tau)$ has uniformly bounded curvature and $M_\tau$ has uniformly bounded second fundamental form.

Define $\hat{ \Sigma } = \lrbrace{ \lrrb{ x, \tau }: x \in M_t, \tau \in (0,T] }$ to be the space-time track of $M_\tau$ within the canonical Ricci shrinker $\lrrb{O^{n+1} \times (0, T], \hat{ g } }$. Then $\hat{ \Sigma }$ is an approximate soliton, i.e.

\begin{equation}
 \nonumber
 \check{ E }_N = H^{ \hat{ \Sigma } }  + \nu^{ \hat{ \Sigma } }   \hat{f} \approx 0
\end{equation}

where $H^{ \hat{ \Sigma } }$ and $\nu^{ \hat{ \Sigma } }$ denote the mean curvatue and the normal vector of $\hat{ \Sigma }$. By this we mean that $N \abs{\hat{ E }_N}_{\hat{ g } }$ is bounded uniformly locally on $\hat{ \Sigma }$ (independently of $N$).

\end{theorem}

The proof of theorem \ref{thm:canonicalmcfsolitonsinrfshrinkers} consists of computing the relevant quantities (see appendix \ref{appendix:canonicalmcfsolitonsexpandingrfbck}). We keep the terminology and call $\hat{ \Sigma }$ a canonical soliton of the mean curvature flow (within a canonical shrinking Ricci soliton).

\subsection{Canonical solitons within steady Ricci flow solitons}

We consider solutions $\lrrb{M^n, M_\tau, O^{n+1}, g^O (\tau)}$ of (\ref{eq:meancurvatureflowinaricciflowbackground}) for $S = 2Ric_{g^O}$ and we choose $I = (0,T]$.

Associated to the backwards Ricci flow $g^O (\tau)$ we have a canonical steady Ricci soliton which can be described as follows (see \cite{estherpeter1}, theorem 5.1). Suppose $\lrrb{O^{n+1}, g^O (\tau)}$, $\tau \in (0,T]$, has uniformly bounded curvature. For a large enough parameter $N > 0$ we define the metric $\bar{g}^O$ on $O^{n+1} \times (0, T]$ by

\begin{equation}
 \label{eq:canonicalriccisteady}
 \bar{g}_{ij} = \begin{cases}
                     g^O_{ij}, \; i,j \geq 0 \\
		     0, \; i \geq 1, j = 0 \\
		     N + R \lrrb{ g^O }
                    \end{cases}
\end{equation}

where $0$ denotes the time direction on $O^{n+1} \times (0, T]$ and $R_{g^O}$ denotes the scalar curvature.

Then up to errors of order $\frac{1}{N}$ the metric $\bar{g}$ is an approximate gradient steady Ricci soliton:

\begin{equation}
 \nonumber
 E_N = Ric_{\bar{g}} + Hess_{\bar{g}} \lrrb{ - N \tau} \approx 0
\end{equation}

by which we mean that the quantity $N \abs{ E_N }_{ \bar{ g } }$ is bounded uniformly locally on $O^{n+1} \times (0, T]$ (independently of $N$).

We denote the potential of the canonical steady Ricci soliton by $\bar{ f } \lrrb{y , \tau} = - N \tau$. As for the expanding and shrinking cases we can show that there is a natural corresponding construction for the coupled flow.

\begin{theorem}
 \label{thm:canonicalmcfsolitonsinrfsteady}

 Suppose $\lrrb{M^n, M_\tau, O^{n+1}, g^O (\tau)}$ is a solution of (\ref{eq:meancurvatureflowinaricciflowbackground}) such that $g^O(\tau)$ has uniformly bounded curvature and $M_\tau$ has uniformly bounded second fundamental form.

Define $\bar{ \Sigma } = \lrbrace{ \lrrb{ x, \tau }: x \in M_\tau, \tau \in (0,T] }$ to be the space-time track of $M_\tau$ within the canonical steady Ricci soliton $\lrrb{O^{n+1} \times (0, T], \bar{g}}$. Then $\bar{ \Sigma }$ is an approximate soliton, i.e.

\begin{equation}
 \nonumber
 \check{ E }_N = H^{ \bar{ \Sigma } }  + \nu^{ \bar{ \Sigma } }   \bar{f} \approx 0
\end{equation}

where $H^{ \bar{ \Sigma } }$ and $\nu^{ \bar{ \Sigma } }$ denote the mean curvatue and the normal vector of $\bar{ \Sigma }$. By this we mean that $N \abs{\bar{ E }_N}_{ \bar{ g } }$ is bounded uniformly locally on $\bar{ \Sigma }$ (independently of $N$).

\end{theorem}

The proof of theorem \ref{thm:canonicalmcfsolitonsinrfsteady} consists of computing the relevant quantities (see appendix \ref{appendix:canonicalmcfsolitonsexpandingrfbck}). We keep the terminology and call $\bar{ \Sigma }$ a canonical soliton of the mean curvature flow (within a canonical steady Ricci soliton).

\section{Canonical solitons and Harnack quantities}

An application of canonical solitons for geometric flows is to generate Harnack quantities. The general idea is that a differential Harnack estimate corresponds to a preserved curvature condition on a canonical soliton.

E. Cabezas-Rivas and P. Topping have shown in \cite{estherpeter2} that this procedure works for the Ricci flow. The Ricci curvature of the canonical expanding Ricci soliton converges as $N \rightarrow \infty$ to Hamilton's matrix Harnack quantity for the Ricci flow (see \cite{hamiltonrf}). Moreover preserved curvature conditions on the canonical expanding Ricci soliton recover S. Brendle's Harnack quantity (see \cite{brendleharnack}) and also provide new Harnack quantities. These correspond to the curvature operator of the canonical soliton lying in certain cones.

It has been shown in \cite{kotschwarharnack} that the same approach works also for curvature flows in Euclidean space. The corresponding Harnack quantities which are due to R. Hamilton, B. Andrews and K. Smoczyk (see \cite{hamiltonmcf}, \cite{andrewsharnack} and \cite{smoczykharnack}) can be recovered as the second fundamental form of a canonical self-expander.

Our construction of canonical mean curvature flow solitons within expanding Ricci flow solitons from theorem \ref{thm:canonicalmcfsolitonsinrfexpanders} provides a link between the Harnack estimate for the mean curvature flow and the Harnack estimate for the Ricci flow. From (\ref{eq:canonicalmcfsolitoninrfexpanderssecondffuptoerrors}) we get

\begin{proposition}
 \label{prop:canonicalmcfsolitoninrfbackgroundsecondfflimit}
 The second fundamental from of the canonical soliton $\check{ \Sigma }$ converges as $N \rightarrow \infty$ to a limit $h^{ \check{ \Sigma } }_\infty$ with

\begin{eqnarray}
 \label{eq:canonicalmcfsolitoninrfbackgroundsecondfflimit}
 h^{ \check{ \Sigma } }_\infty \lrrb{V + \frac{\partial}{\partial t}, V + \frac{\partial}{\partial t}} = \\
 \nonumber
  \frac{\partial}{\partial t} H + h \lrrb{V, V} + \frac{H}{2t} + 2 \skp{ V, \nabla H } + \\
 \nonumber
  2 Ric_{ g^O } \lrrb{ V, \nu } - H Ric_{ g^O } \lrrb{\nu, \nu} + \frac{1}{2} \check{ \nabla }_{ \nu } R_{ g^O }.
\end{eqnarray}

for any vector field $V$ on $M^n$.

\end{proposition}

On the other hand we have (see \cite[proposition 2.6]{estherpeter2})

\begin{proposition}
 \label{prop:canonicalricciexpanderriccilimit}
 The Ricci curvature of the canonical expanding Ricci soliton $\check{g}$ converges as $N \rightarrow \infty$ to $\lrrb{Ric_{\check{g}}}_\infty$ with

 \begin{eqnarray}
  \label{eq:canonicalricciexpanderriccilimit}
  \lrrb{Ric_{\check{g}}}_\infty \lrrb{ X + \frac{\partial}{\partial t}, X + \frac{\partial}{\partial t} } = \\
  \nonumber
  Ric_{g^O} \lrrb{X, X} + g^O \lrrb{ X, {}^{g^O} \nabla R_{g^O} } + \frac{1}{2} \lrrb{ \frac{ \partial R_{g^O} }{ \partial t } + \frac{ R_{g^O} }{ t } }
 \end{eqnarray}

for any vector field $X$ on $O^{n+1}$.

\end{proposition}

From proposition \ref{prop:canonicalricciexpanderriccilimit} we see that the Ricci curvature of the canonical expanding Ricci soliton converges to R. Hamilton's Harnack quantity for the Ricci flow (see \cite{hamiltonrf})

\begin{equation}
  \label{eq:rfharnack}
  Z \lrrb{ X, X } = Ric_{g^O} \lrrb{X, X} + g^O \lrrb{ X, {}^{g^O} \nabla R_{g^O} } + \frac{1}{2} \lrrb{ \frac{ \partial R_{g^O} }{ \partial t } + \frac{ R_{g^O} }{ t } }.
\end{equation}

 as $N \rightarrow \infty$. If $O^{n+1} = \set{R}^{n+1}$ with the Euclidean metric, then we see from proposition \ref{prop:canonicalmcfsolitoninrfbackgroundsecondfflimit} that the second fundamental form of the mean curvature flow soliton within the canonical expanding Ricci soliton converges as $N \rightarrow \infty$ to Hamilton's Harnack quantity for the mean curvature flow

 \begin{equation}
   \label{eq:mcfharnack}
   \tilde{Z} \lrrb{ V, V } =  \frac{\partial}{\partial t} H + h \lrrb{V, V} + 2 \skp{ V, \nabla H } + \frac{H}{2t}
 \end{equation}

 (see \cite{hamiltonmcf}), providing a link between the Harnack quantity for the Ricci flow and the Harnack quantity for the mean curvature flow.

This seems to suggest that the quantity (\ref{eq:canonicalmcfsolitoninrfbackgroundsecondfflimit}) should have a similar meaning for the coupled flow. Unfortunately no suitable preserved curvature condition for proving a Harnack type estimate for the coupled flow is known

\section{Canonical solitons and the modified $\mathcal{F}$-functional}

In view of the relation of canonical solitons to Harnack estimates it seems natural that the curvature of canonical solitons also plays an important role for the coupled flow (\ref{eq:definitionmeancurvatureflowinaricciflowbackground}). Again we focus here on the case $S = - 2Ric$, so we get canonical expanding Ricci solitons and (within these) canonical mean curvature flow solitons. 

Nonnegative second fundamental form is not preserved in general for the coupled flow, even for a Ricci soliton background. Hence we can not expect a direct analogue of the differential Harnack estimate for the mean curvature flow. However we want to point out a link between canonical solitons and a version of the $\mathcal{F}$-functional (see \cite[section 1.1]{perelman1}).

In \cite{eckerentropyharnack} K. Ecker found a remarkable connection between the $\mathcal{W}$-functional for the Ricci flow and the differential Harnack expression for the mean curvature flow in $\set{R}^{n+1}$. Suppose that for a solution $\lrrb{M_t}_{t \in [0,T)}$ of the mean curvature flow in $\set{ R }^{n+1}$ the hypersurfaces $M_t$ are the boundary of open bounded sets $\Omega_t$ in ${ R }^{n+1}$. One can consider the integral over $\Omega_t$ of G. Perelman's $\mathcal{W}$-integrand (see \cite[proposition 9.1]{perelman1})

\begin{equation}
  \nonumber
  \mathcal{ W } \lrrb{ \Omega_t, f, \tau {\lrrb{ t } } } = \int_{\Omega_t } \lrrb{ \tau \abs{ \nabla f }^2 + f - (n+1) } u dV + 2 \tau \int_{M_t} H u dA
\end{equation}

\footnote{$dV$ denotes the volume density in $\set{ R }^{n+1}$ and $dA$ the associated area density of a boundary}  for a reverse time parameter $\frac{ d }{ dt } \tau \lrrb{ t } = -1$, where the $\mathcal{W}$-integrand is itself defined using a positive solution of the backwards heat equation

\begin{eqnarray}
   \nonumber
   u = \frac{ e^{-f} }{  \lrrb{ 4 \pi \tau }^{ \frac{ n+1 }{ 2 } } } \\
   \frac{ \partial }{ \partial t } f + \triangle f = \abs{ \nabla f }^2 + \frac{ n + 1}{ 2 \tau }
\end{eqnarray}

in $\Omega_t$ with von Neumann boundary condition

\begin{equation}
  \nonumber
  H = \skp{ \nabla f, \nu }
\end{equation}

on $M_t = \partial \Omega_t$.

The boundary term of the evolution of $\mathcal{ W }$ then consists of Hamilton's Harnack quantity for the mean curvature flow (\ref{eq:mcfharnack}):

\begin{equation}
  \label{eq:eckerwevolution}
  \frac{ d }{ dt } \mathcal{ W } \lrrb{ \Omega_t, f , \tau {\lrrb{ t } } } = 2\tau \int_{\Omega_t} \abs{ \nabla_i \nabla_j f - \frac{ \delta_i^j }{ 2 \tau } }^2 dV + 2 \tau \int_{M_t} \tilde{ Z } \lrrb{ -\nabla f, -\nabla f } u dA
\end{equation}

(see \cite[proposition 3.4]{eckerentropyharnack}).

Furthermore K. Ecker conjectured that the boundary term in (\ref{eq:eckerwevolution}) is nonnegative for any compact solution $M_t = \partial \Omega_t$ of the mean curvature flow in $\set{R}^{n+1}$ with positive mean curvature. This conjecture is still open.

J. Lott showed in \cite{lottmcfrfbckgr} that analogous $\mathcal{ F }$-versions of these relations hold for solutions of (\ref{eq:definitionmeancurvatureflowinaricciflowbackground}). He defined a version $I_\infty$ of G. Perelman's $\mathcal{F}$-functional for manifolds $O^{n+1}$ with boundary by adding an appropriate boundary term \footnote{we mean here appropriate from a variational viewpoint} to the original integral. For $f \in C^\infty \lrrb{O^{n+1}}$ we have

\begin{equation}
 \label{eq:ghyfunctional}
 I_\infty \lrrb{g^O, f} = \int_{O^{n+1}} R_{g^O}^\infty e^{-f} dV_{g^O} + 2 \int_{\partial O^{n+1}} H^\infty e^{-f} dA_{g^O}
\end{equation}

where

\begin{equation}
  \nonumber
  R_{g^O}^\infty = R_{g^O} + 2 \triangle_{g^O} f - \abs{ {}^{g^O} \nabla f }^2_{g^O}
\end{equation}

is the analogue of the scalar curvature on the smooth metric measure space $\mathcal{O} = \lrrb{ O^{n+1}, g^O, e^{-f} dV_{g^O} }$ and

\begin{equation}
  \nonumber
  H^\infty = H - \nu f
\end{equation}

is the analogue of the mean curvature. Here $\nu$ denotes the outward-pointing unit normal on $\partial O^{n+1}$ and $dA_{g^O}$ the induced area on $\partial O^{n+1}$. $I_\infty$ can be seen as a weighted \footnote{weighted with respect to the measure $e^{-f} dV_{g^O}$} version of the Gibbons-Hawking-York functional (see \cite{gibbonshawkingfunc} and \cite{yorkfunc})

\begin{equation}
  \nonumber
  I_{GHY} \lrrb{ g^O } = \int_{O^{n+1}} R_{g^O} e^{-f} dV_{g^O} + 2 \int_{ \partial O^{n+1} } H dA_{g^O}
\end{equation}

which occurs in quantum gravity. 

Now if $g^O(t)$ is a solution of the modified Ricci flow

\begin{equation}
  \nonumber
  \frac{\partial}{\partial t} g^O = -2 \lrrb{ Ric_{g^O} + Hess_{g^O} \lrrb{f} }
\end{equation}

for a solution $u = e^{-f}$ of the conjugate heat equation

\begin{equation}
  \nonumber
  \frac{\partial}{\partial t} u = \lrrb{- \triangle_{g^O} + R_{g^O}} u
\end{equation}

with boundary condition

\begin{equation}
  \nonumber
  H - \nu f = 0
\end{equation}

then we have (see \cite[theorem 1]{lottmcfrfbckgr})

\begin{eqnarray}
 \label{eq:ghyevolutionequation}
 \frac{d}{dt} I_\infty = 2 \int_{O^{n+1}} \abs{ Ric_{g^O} + Hess_{g^O} \lrrb{f} }^2 e^{-f} dV_{g^O} + \\
 \nonumber
 2 \int_{\partial O^{n+1}} \left( \frac{\partial}{\partial t} H - 2 \skp{ \nabla f, \nabla H } + h \lrrb{ \nabla f, \nabla f } \right. \\
 \nonumber
 \left. - 2 Ric_{g^O} \lrrb{ \nu, \nabla f } + \frac{1}{2} \nu R_{g^O} - H Ric_{g^O} \lrrb{ \nu, \nu } \right) e^{-f} dA_{g^O}
\end{eqnarray}

where $\nabla$ denotes the derivative on $\partial O^{n+1}$ and $h$, $H$ denote the second fundamental form and the mean curvature of $\partial O^{n+1}$.

From proposition \ref{prop:canonicalmcfsolitoninrfbackgroundsecondfflimit} we see that up to the $\frac{H}{2t}$ term the limit second fundamental form of our canonical soliton evaluated at $- \nabla f$ exactly matches the boundary integrand of (\ref{eq:ghyevolutionequation}). The former is $0$ on mean curvature flow solitons in gradient expanding Ricci solitons, whereas the latter is $0$ on mean curvature flow solitons in gradient steady Ricci solitons (see \cite[proposition 7]{lottmcfrfbckgr}).

In the case of Euclidean ambient space $\set{R}^{n+1}$ (\ref{eq:ghyevolutionequation}) is the $\mathcal{F}$-version of K. Ecker's result (see \cite[proposition 3.2, 3.4]{eckerentropyharnack}). In this case the limit second fundamental form of the canonical soliton from theorem \ref{thm:canonicalmcfsolitonsinrfexpanders} exactly matches Hamilton's Harnack quantity for the mean curvature flow (\ref{eq:mcfharnack}), which in turn matches the boundary integrand in (\ref{eq:ghyevolutionequation}) up to the $\frac{H}{2t}$ term.

In the same fashion as proposition \ref{prop:canonicalmcfsolitoninrfbackgroundsecondfflimit} the canonical solitons from theorem \ref{thm:canonicalmcfsolitonsinrfshrinkers} and theorem \ref{thm:canonicalmcfsolitonsinrfsteady}  yield analogous quantities to the boundary integrand in (\ref{eq:ghyevolutionequation}) for a backward Ricci flow background.

\appendix

\section{Computations 1: canonical mean curvature flow solitons within canonical expanding Ricci flow solitons}

\label{appendix:canonicalmcfsolitonsexpandingrfbck}

We prove theorem \ref{thm:canonicalmcfsolitonsinrfexpanders} by computing the induced quantities on $\check{ \Sigma }$ which we denote by $\cdot ^{ \check{ \Sigma } }$. Points on $\check{ \Sigma }$ are denoted by $z = F ^{ \check{ \Sigma } } \lrrb{ x, t} = \lrrb{ F_t (x), t}$. In index notation $0$ always denotes the time direction and we also use $x^0 = t$ here. As before we use bare letters to denote quantities on the hypersurfaces $M_t$, e.g. $g$ or $h$. By $A \approx B$ up to errors of order $\frac{ 1 }{ N }$ on $\check{ \Sigma }$ we always mean $A = B + E_N$ and $N \abs{ E_N }_{ \check{ g } }$ is bounded uniformly locally on $\check{ \Sigma }$ (independently of $N$).

The computations of the Christoffel symbols and the curvature of $\check{ g }$ can be found in \cite{estherpeter2}. To compute the second fundamental form of $\check{ \Sigma }$ we need the Christoffel symbols of $\check{ g }$, which are given by:

\begin{eqnarray}
 \label{eq:canonicalmcfsolitoninrfexpanderschristoffel}
 \check{ \Gamma }^\alpha_{\beta \gamma} = \lrrb{ \Gamma_{g^O} }^\alpha_{\beta \gamma}, \;
 \check{ \Gamma }^\alpha_{\beta 0} = - \lrrb{ \lrrb{Ric_{g_O}}^\alpha_\beta + \frac{\delta^\alpha_\beta}{2t} }, \; 
 \check{ \Gamma }^\alpha_{0 0} = - \frac{1}{2} \lrrb{g^O}_{\alpha \beta} \frac{\partial R_{g^O}}{\partial y^\beta} \\
 \nonumber
 \check{ \Gamma }^0_{\beta \gamma} = \check{ g }_{00}^{-1} \lrrb{ \frac{ \lrrb{Ric_{g^O}}_{ \beta \gamma } }{t} + \frac{ g^O_{\beta \gamma }}{2 t^2} } , \;
 \check{ \Gamma }^0_{\beta 0} = \frac{1}{2t} \check{ g }_{00}^{-1} \frac{\partial R_{g^O}}{\partial y^\beta} \\
 \nonumber
 \check{ \Gamma }^0_{0 0} = - \frac{3}{2t} + \frac{\check{ g }_{00}^{-1}}{2t} \lrrb{ \frac{ R_{g^O}}{t} + \frac{\partial R_{g^O}}{\partial t} + \frac{n + 1}{2 t^2} }.
\end{eqnarray}

A basis for $T_z \check{ \Sigma }$ is given by

\begin{eqnarray}
 \nonumber
 \lrbrace{ \lrrb{F_t}_* \frac{\partial}{\partial x^1}, \ldots, \lrrb{F_t}_* \frac{\partial}{\partial x^n}, - H \nu + \frac{\partial}{\partial t} } = \\
 \nonumber \lrbrace{ \frac{\partial}{\partial x^1} F_t, \ldots, \frac{\partial}{\partial x^n} F_t, - H \nu + \frac{\partial}{\partial t} }
\end{eqnarray}

Therefore we get

\begin{equation}
 \label{eq:canonicalmcfsolitoninrfexpandersnormal}
 \nu ^{ \check{ \Sigma } } = \frac{1}{\sigma_N} \lrrb{\nu + \frac{H}{ \check{ g }_{00} t } \frac{\partial}{\partial t}}
\end{equation}

where $\sigma_N = \sqrt{ t^{-1} + H^2 t^{-2} \check{ g }_{00}^{-1} }$.

(\ref{eq:canonicalmcfsolitoninrfexpandersnormal}) yields

\begin{equation}
 \label{eq:canonicalmcfsolitoninrfexpanderspotentialgradient}
 \nu ^{ \check{ \Sigma } } \check{ f} = \frac{HN}{2 t^3 \check{ g }_{00} \sigma_N}.
\end{equation}

From 

\begin{equation}
 \nonumber
 g ^{ \check{ \Sigma } }_{ij} = \check{ g } \lrrb{ \frac{\partial}{\partial x^i} F ^{\check{ \Sigma }}, \frac{\partial}{\partial x^j} F ^{\check{ \Sigma }} }
\end{equation}

we obtain

\begin{equation}
 \label{eq:canonicalmcfsolitoninrfexpandersmetric}
 g ^{ \check{ \Sigma } }_{ij} = \begin{cases}
                                 \frac{1}{t} g_{ij}, \; i,j \geq 0\\
				 0, \; i \geq 0, j = 0 \\
				 \frac{H^2}{t} + \check{g}_{00}, \; i = j = 0.
                                \end{cases}
\end{equation}

For the inverse metric we have then

\begin{equation}
 \label{eq:canonicalmcfsolitoninrfexpandersinversemetric}
 \lrrb{ g ^{ \check{ \Sigma } } }^{ij} = \begin{cases}
                                 t g^{ij}, \; i,j \geq 0\\
				 0, \; i \geq 0, j = 0 \\
				 \frac{ t }{H^2 + t \check{g}_{00}}, \; i = j = 0.
                                \end{cases}
\end{equation}

Hence the inverse metric $\lrrb{ g ^{ \check{ \Sigma } } }^{-1}$ converges as $N \rightarrow \infty$ to the degenerate metric $\lrrb{ g ^{ \check{ \Sigma } }_\infty }^{-1}$ given by

\begin{equation}
 \label{eq:canonicalmcfsolitoninrfexpandersinverselimitmetric}
 \lrrb{ g ^{ \check{ \Sigma } }_\infty }^{ij} = \begin{cases}
                                 t g^{ij}, \; i,j \geq 0\\
				 0, \; i \geq 0, j = 0 \\
				 0, \; i = j = 0.
                                \end{cases}
\end{equation}

Using the formula

\begin{equation}
 \nonumber
 h^{ \check{ \Sigma } }_{ij} = - \check{ g } \lrrb{ \frac{\partial}{\partial x^i} F ^{\check{ \Sigma }}, \frac{\partial}{\partial x^j} F ^{\check{ \Sigma }} }
\end{equation}

we can compute the second fundamental form of $\check{ \Sigma }$ (see also (\ref{eq:canonicalmcfsolitoninrfexpanderschristoffel}))

\begin{equation}
 \label{eq:canonicalmcfsolitoninrfexpanderssecondff}
 h^{ \check{ \Sigma } }_{ij} = \frac{1}{t\sigma_N} \begin{cases}
                                 h_{ij} + \frac{ H }{  t\check{ g }_{00} } \lrrb{ Ric_{ g^O } \lrrb{ \frac{\partial}{\partial x^i} F_t, \frac{\partial}{\partial x^j} F_t} + \frac{1}{2t} g_{ij} }, \; i,j \geq 1 \\
				 \frac{\partial}{\partial x^i} H + Ric_{ g^O } \lrrb{ \frac{\partial}{\partial x^i} F_t, \nu } - \frac{H}{2t \check{g}_{00}} \frac{\partial}{\partial x^i} F_t \lrrb{ R_{ g^O } }, \; i \geq 1, j = 0 \\
				 \frac{\partial}{\partial t} H + \frac{H}{2t} - H Ric_{ g^O } \lrrb{\nu, \nu} + \frac{1}{2} \check{ \nabla }_{ \nu } R_{ g^O } + \\ 
				  \frac{H}{t \check{g}_{00}} \left ( \frac{H^2}{2t} + H^2 Ric_{ g^O } \lrrb{\nu, \nu} - H^2 \check{ \nabla }_{ \nu } R_{ g^O } - \right . \\
				  \left . \frac{1}{t} R_{ g^O } - \frac{1}{2} \frac{\partial}{\partial t} R_{ g^O } + \frac{n}{4t^2} \right ), \; i = j = 0.
                               \end{cases}
\end{equation}

Up to errors of order $\frac{1}{N}$ this means

\begin{equation}
 \label{eq:canonicalmcfsolitoninrfexpanderssecondffuptoerrors}
 h^{ \check{ \Sigma } }_{ij} \approx \frac{1}{t\sigma_N} \begin{cases}
                                 h_{ij}, \; i,j \geq 1 \\
				 \frac{\partial}{\partial x^i} H + Ric_{ g^O } \lrrb{ \frac{\partial}{\partial x^i} F_t, \nu }, \; i \geq 1, j = 0 \\
				 \frac{\partial}{\partial t} H + \frac{H}{2t} - H Ric_{ g^O } \lrrb{\nu, \nu} + \frac{1}{2} \check{ \nabla }_{ \nu } R_{ g^O }, \; i = j = 0.
                               \end{cases}
\end{equation}

Hence (\ref{eq:canonicalmcfsolitoninrfexpandersinversemetric}) and (\ref{eq:canonicalmcfsolitoninrfexpanderssecondffuptoerrors}) show that the mean curvature of $\check{ \Sigma }$ is given by (up to errors of order $\frac{1}{N}$)

\begin{equation}
 \label{eq:canonicalmcfsolitoninrfexpandersmeancurvature}
 H^{ \check{ \Sigma } } \approx \frac{t}{\sigma_N} H \approx \sqrt{t} H.
\end{equation}

Together (\ref{eq:canonicalmcfsolitoninrfexpanderspotentialgradient}) and (\ref{eq:canonicalmcfsolitoninrfexpandersmeancurvature}) then yield

\begin{equation}
 \label{eq:canonicalmcfsolitoninrfexpanderssolitonequation}
 H^{ \check{ \Sigma } } \approx \nu^{ \check{ \Sigma } } \check{ f }
\end{equation}

from which we can deduce theorem \ref{thm:canonicalmcfsolitonsinrfexpanders}.

\section{Computations 2: canonical mean curvature flow solitons within canonical shrinking Ricci flow solitons}

\label{appendix:canonicalmcfsolitonsshrinkingrfbck}

We use analogue assumptions and notations as in appendix \ref{appendix:canonicalmcfsolitonsexpandingrfbck} to prove theorem \ref{thm:canonicalmcfsolitonsinrfshrinkers}.

The Christoffel symbols of $\hat{ g }$ are given by (see \cite{estherpeter1}).

\begin{eqnarray}
 \label{eq:canonicalmcfsolitoninrfshrinkerschristoffel}
 \hat{ \Gamma }^\alpha_{\beta \gamma} = \lrrb{ \Gamma_{g^O} }^\alpha_{\beta \gamma}, \;
 \hat{ \Gamma }^\alpha_{\beta 0} = \lrrb{Ric_{g_O}}^\alpha_\beta - \frac{\delta^\alpha_\beta}{2\tau}, \; 
 \hat{ \Gamma }^\alpha_{0 0} = - \frac{1}{2} \lrrb{g^O}_{\alpha \beta} \frac{\partial R_{g^O}}{\partial y^\beta} \\
 \nonumber
 \hat{ \Gamma }^0_{\beta \gamma} = \hat{ g }_{00}^{-1} \lrrb{ - \frac{ \frac{ g^O_{\beta \gamma }}{2 \tau^2} - \lrrb{Ric_{g^O}}_{ \beta \gamma } }{\tau} } , \;
 \hat{ \Gamma }^0_{\beta 0} = \frac{1}{2\tau} \hat{ g }_{00}^{-1} \frac{\partial R_{g^O}}{\partial y^\beta} \\
 \nonumber
 \hat{ \Gamma }^0_{0 0} = - \frac{3}{2\tau} + \frac{\hat{ g }_{00}^{-1}}{2\tau} \lrrb{ \frac{ R_{g^O}}{\tau} + \frac{\partial R_{g^O}}{\partial \tau} + \frac{n + 1}{2 \tau^2} }.
\end{eqnarray}

We can proceed as in appendix \ref{appendix:canonicalmcfsolitonsexpandingrfbck} to compute the relevant quantities.

For the normal vector we get

\begin{equation}
 \label{eq:canonicalmcfsolitoninrfshrinkersnormal}
 \nu ^{ \hat{ \Sigma } } = \frac{1}{\sigma_N} \lrrb{\nu + \frac{H}{ \hat{ g }_{00} \tau } \frac{\partial}{\partial \tau}}
\end{equation}

where $\sigma_N = \sqrt{ \tau^{-1} + H^2 \tau^{-2} \hat{ g }_{00}^{-1} }$.

(\ref{eq:canonicalmcfsolitoninrfshrinkersnormal}) yields

\begin{equation}
 \label{eq:canonicalmcfsolitoninrfshrinkerspotentialgradient}
 \nu ^{ \hat{ \Sigma } } \hat{ f} = - \frac{HN}{2 \tau^3 \hat{ g }_{00} \sigma_N}.
\end{equation}

The induced metric on $\hat{ \Sigma }$ is

\begin{equation}
 \label{eq:canonicalmcfsolitoninrfshrinkersmetric}
 g ^{ \hat{ \Sigma } }_{ij} = \begin{cases}
                                 \frac{1}{\tau} g_{ij}, \; i,j \geq 0\\
				 0, \; i \geq 0, j = 0 \\
				 \frac{H^2}{\tau} + \hat{g}_{00}, \; i = j = 0.
                                \end{cases}
\end{equation}

For the inverse metric we have then

\begin{equation}
 \label{eq:canonicalmcfsolitoninrfshrinkersinversemetric}
 \lrrb{ g ^{ \hat{ \Sigma } } }^{ij} = \begin{cases}
                                 \tau g^{ij}, \; i,j \geq 0\\
				 0, \; i \geq 0, j = 0 \\
				 \frac{ \tau }{H^2 + \tau \hat{g}_{00}}, \; i = j = 0.
                                \end{cases}
\end{equation}

Hence the inverse metric $\lrrb{ g ^{ \hat{ \Sigma } } }^{-1}$ converges as $N \rightarrow \infty$ to the degenerate metric $\lrrb{ g ^{ \hat{ \Sigma } }_\infty }^{-1}$ given by

\begin{equation}
 \label{eq:canonicalmcfsolitoninrfshrinkersinverselimitmetric}
 \lrrb{ g ^{ \hat{ \Sigma } }_\infty }^{ij} = \begin{cases}
                                 \tau g^{ij}, \; i,j \geq 0\\
				 0, \; i \geq 0, j = 0 \\
				 0, \; i = j = 0.
                                \end{cases}
\end{equation}

We can compute the second fundamental form of $\hat{ \Sigma }$ (see also (\ref{eq:canonicalmcfsolitoninrfshrinkerschristoffel}))

\begin{equation}
 \label{eq:canonicalmcfsolitoninrfshrinkerssecondff}
 h^{ \hat{ \Sigma } }_{ij} = \frac{1}{\tau \sigma_N} \begin{cases}
                                 h_{ij} + \frac{ H }{  \tau\hat{ g }_{00} } \lrrb{ - Ric_{ g^O } \lrrb{ \frac{\partial}{\partial x^i} F_\tau, \frac{\partial}{\partial x^j} F_\tau} + \frac{1}{2\tau} g_{ij} }, \; i,j \geq 1 \\
				 \frac{\partial}{\partial x^i} H - Ric_{ g^O } \lrrb{ \frac{\partial}{\partial x^i} F_\tau, \nu } - \frac{H}{2\tau \hat{g}_{00}} \frac{\partial}{\partial x^i} F_\tau \lrrb{ R_{ g^O } }, \; i \geq 1, j = 0 \\
				 \frac{\partial}{\partial \tau} H + \frac{H}{2\tau} + H Ric_{ g^O } \lrrb{\nu, \nu} + \frac{1}{2} \hat{ \nabla }_{ \nu } R_{ g^O } + \\ 
				  \frac{H}{\tau \hat{g}_{00}} \left ( - \frac{H^2}{2\tau} + H^2 Ric_{ g^O } \lrrb{\nu, \nu} - H^2 \hat{ \nabla }_{ \nu } R_{ g^O } - \right . \\
				  \left . \frac{1}{\tau} R_{ g^O } - \frac{1}{2} \frac{\partial}{\partial \tau} R_{ g^O } + \frac{n}{4\tau^2} \right ), \; i = j = 0.
                               \end{cases}
\end{equation}

Up to errors of order $\frac{1}{N}$ this means

\begin{equation}
 \label{eq:canonicalmcfsolitoninrfshrinkerssecondffuptoerrors}
 h^{ \hat{ \Sigma } }_{ij} \approx \frac{1}{\tau\sigma_N} \begin{cases}
                                 h_{ij}, \; i,j \geq 1 \\
				 \frac{\partial}{\partial x^i} H - Ric_{ g^O } \lrrb{ \frac{\partial}{\partial x^i} F_\tau, \nu }, \; i \geq 1, j = 0 \\
				 \frac{\partial}{\partial \tau} H + \frac{H}{2\tau} + H Ric_{ g^O } \lrrb{\nu, \nu} + \frac{1}{2} \hat{ \nabla }_{ \nu } R_{ g^O }, \; i = j = 0.
                               \end{cases}
\end{equation}

Hence (\ref{eq:canonicalmcfsolitoninrfshrinkersinversemetric}) and (\ref{eq:canonicalmcfsolitoninrfshrinkerssecondffuptoerrors}) show that the mean curvature of $\hat{ \Sigma }$ is given by (up to errors of order $\frac{1}{N}$)

\begin{equation}
 \label{eq:canonicalmcfsolitoninrfshrinkersmeancurvature}
 H^{ \hat{ \Sigma } } \approx \frac{\tau}{\sigma_N} H \approx \sqrt{\tau} H.
\end{equation}

Together (\ref{eq:canonicalmcfsolitoninrfshrinkerspotentialgradient}) and (\ref{eq:canonicalmcfsolitoninrfshrinkersmeancurvature}) then yield

\begin{equation}
 \label{eq:canonicalmcfsolitoninrfshrinkerssolitonequation}
 H^{ \hat{ \Sigma } } \approx - \nu^{ \hat{ \Sigma } } \hat{ f }
\end{equation}

from which we can deduce theorem \ref{thm:canonicalmcfsolitonsinrfshrinkers}.

\section{Computations 3: canonical mean curvature flow solitons within canonical steady Ricci flow solitons}

\label{appendix:canonicalmcfsolitonssteadyrfbck}

We use analogue assumptions and notations as in appendix \ref{appendix:canonicalmcfsolitonsexpandingrfbck} to prove theorem \ref{thm:canonicalmcfsolitonsinrfsteady}.

The Christoffel symbols of $\bar{ g }$ are given by (see \cite{estherpeter1}).

\begin{eqnarray}
 \label{eq:canonicalmcfsolitoninrfsteadychristoffel}
 \bar{ \Gamma }^\alpha_{\beta \gamma} = \lrrb{ \Gamma_{g^O} }^\alpha_{\beta \gamma}, \;
 \bar{ \Gamma }^\alpha_{\beta 0} = \lrrb{Ric_{g_O}}^\alpha_\beta, \; 
 \bar{ \Gamma }^\alpha_{0 0} = - \frac{1}{2} \lrrb{g^O}_{\alpha \beta} \frac{\partial R_{g^O}}{\partial y^\beta} \\
 \nonumber
 \bar{ \Gamma }^0_{\beta \gamma} = - \frac{ 1 }{ N + R_{g^O} } \lrrb{ Ric_{g^O} }_{ \beta \gamma } , \;
 \bar{ \Gamma }^0_{\beta 0} = \frac{1}{2} \frac{\partial R_{g^O}}{\partial y^\beta}, \;
 \bar{ \Gamma }^0_{0 0} = \frac{1}{2} \frac{\partial R_{g^O}}{\partial \tau}.
\end{eqnarray}

We can proceed as in appendix \ref{appendix:canonicalmcfsolitonsexpandingrfbck} to compute the relevant quantities.

For the normal vector we get

\begin{equation}
 \label{eq:canonicalmcfsolitoninrfsteadynormal}
 \nu ^{ \bar{ \Sigma } } = \frac{1}{\sigma_N} \lrrb{\nu + \frac{H}{ \bar{ g }_{00} } \frac{\partial}{\partial \tau}}
\end{equation}

where $\sigma_N = \sqrt{ 1 + H^2 \bar{ g }_{00}^{-1} }$.

(\ref{eq:canonicalmcfsolitoninrfsteadynormal}) yields

\begin{equation}
 \label{eq:canonicalmcfsolitoninrfsteadypotentialgradient}
 \nu ^{ \bar{ \Sigma } } \bar{ f} = - \frac{HN}{ \bar{ g }_{00} \sigma_N}.
\end{equation}

The induced metric on $\bar{ \Sigma }$ is

\begin{equation}
 \label{eq:canonicalmcfsolitoninrfsteadymetric}
 g ^{ \bar{ \Sigma } }_{ij} = \begin{cases}
                                 g_{ij}, \; i,j \geq 0\\
				 0, \; i \geq 0, j = 0 \\
				 \frac{H^2}+ \bar{g}_{00}, \; i = j = 0.
                                \end{cases}
\end{equation}

For the inverse metric we have then

\begin{equation}
 \label{eq:canonicalmcfsolitoninrfsteadyinversemetric}
 \lrrb{ g ^{ \bar{ \Sigma } } }^{ij} = \begin{cases}
                                 g^{ij}, \; i,j \geq 0\\
				 0, \; i \geq 0, j = 0 \\
				 \frac{ 1 }{H^2 + \tau \bar{g}_{00}}, \; i = j = 0.
                                \end{cases}
\end{equation}

Hence the inverse metric $\lrrb{ g ^{ \bar{ \Sigma } } }^{-1}$ converges as $N \rightarrow \infty$ to the degenerate metric $\lrrb{ g ^{ \bar{ \Sigma } }_\infty }^{-1}$ given by

\begin{equation}
 \label{eq:canonicalmcfsolitoninrfsteadyinverselimitmetric}
 \lrrb{ g ^{ \bar{ \Sigma } }_\infty }^{ij} = \begin{cases}
                                 g^{ij}, \; i,j \geq 0\\
				 0, \; i \geq 0, j = 0 \\
				 0, \; i = j = 0.
                                \end{cases}
\end{equation}

We can compute the second fundamental form of $\bar{ \Sigma }$ (see also (\ref{eq:canonicalmcfsolitoninrfsteadychristoffel}))

\begin{equation}
 \label{eq:canonicalmcfsolitoninrfsteadysecondff}
 h^{ \bar{ \Sigma } }_{ij} = \frac{1}{\sigma_N} \begin{cases}
                                 h_{ij} + \frac{ H }{ N + R_{g^O} } Ric_{ g^O } \lrrb{ \frac{\partial}{\partial x^i} F_\tau, \frac{\partial}{\partial x^j} F_\tau}, \; i,j \geq 1 \\
				 \frac{\partial}{\partial x^i} H - Ric_{ g^O } \lrrb{ \frac{\partial}{\partial x^i} F_\tau, \nu } + \\ \frac{H}{2} \frac{\partial}{\partial x^i} \ln \lrrb{ N + R_{g^O} } + \frac{H^2}{ N + R_{g^O} } Ric_{ g^O } \lrrb{ \frac{\partial}{\partial x^i} F_\tau, \nu }, \; i \geq 1, j = 0 \\
				 \frac{\partial}{\partial \tau} H + H Ric_{ g^O } \lrrb{\nu, \nu} + \\ 
				  \frac{1}{2} \bar{ \nabla }_{ \nu } R_{ g^O } + \frac{H^2}{2} \bar{ \nabla }_{ \nu } \ln \lrrb{ N + R_{g^O} } - \\
				  \frac{H}{2} \frac{\partial}{\partial t} \ln \lrrb{ N + R_{g^O} } - \frac{ H^3 }{ N + R_{g^O} } Ric_{ g^O } \lrrb{\nu, \nu}, \; i = j = 0.
                               \end{cases}
\end{equation}

Up to errors of order $\frac{1}{N}$ this means

\begin{equation}
 \label{eq:canonicalmcfsolitoninrfsteadysecondffuptoerrors}
 h^{ \bar{ \Sigma } }_{ij} \approx \frac{1}{\tau\sigma_N} \begin{cases}
                                 h_{ij}, \; i,j \geq 1 \\
				 \frac{\partial}{\partial x^i} H - Ric_{ g^O } \lrrb{ \frac{\partial}{\partial x^i} F_\tau, \nu }, \; i \geq 1, j = 0 \\
				 \frac{\partial}{\partial \tau} H + H Ric_{ g^O } \lrrb{\nu, \nu} + \frac{1}{2} \bar{ \nabla }_{ \nu } R_{ g^O }, \; i = j = 0.
                               \end{cases}
\end{equation}

Hence (\ref{eq:canonicalmcfsolitoninrfsteadyinversemetric}) and (\ref{eq:canonicalmcfsolitoninrfsteadysecondffuptoerrors}) show that the mean curvature of $\bar{ \Sigma }$ is given by (up to errors of order $\frac{1}{N}$)

\begin{equation}
 \label{eq:canonicalmcfsolitoninrfsteadymeancurvature}
 H^{ \bar{ \Sigma } } \approx \frac{H}{\sigma_N} \approx H.
\end{equation}

Together (\ref{eq:canonicalmcfsolitoninrfsteadypotentialgradient}) and (\ref{eq:canonicalmcfsolitoninrfsteadymeancurvature}) then yield

\begin{equation}
 \label{eq:canonicalmcfsolitoninrfsteadysolitonequation}
 H^{ \bar{ \Sigma } } \approx - \nu^{ \bar{ \Sigma } } \bar{ f }
\end{equation}

from which we can deduce theorem \ref{thm:canonicalmcfsolitonsinrfsteady}.

\end{document}